\documentclass{scrartcl}

\usepackage{graphicx}        
\usepackage{lmodern}

\usepackage{stmaryrd}

\usepackage{amssymb,amsmath,colonequals,bm,nicefrac}
\usepackage{tkz-graph}
\usepackage{url}

\usepackage{pgfplots}
\usepackage{siunitx}
\usepackage{booktabs}
\usepackage[capitalize]{cleveref}
\usetikzlibrary{patterns}
\usepackage{wrapfig}

\newtheorem{Problem}{Problem}

\newcommand\iprod[3][]{\bigl(#2,#3\bigr)_{#1}}
\newcommand\prodf{\iprod[\Omega_f]}
\newcommand\prods{\iprod[\Omega_s]}
\newcommand\prodi{\iprod[\Gamma_i]}
\newcommand\prodd{\iprod[\Gamma_f \setminus \Gamma_f^D]}
\newcommand\norm[1]{\lVert#1\rVert}
\newcommand\normL[1]{\norm{#1}_{L^2}}
\newcommand\jump[1]{\llbracket#1\rrbracket}
\newcommand\jumpn[2]{\jump{\partial_n^{#2}\varphi_{#1}}}
\newcommand\Jump[1]{g_F^{#1}(\varphi_1, \varphi_2)}
\newcommand\bd{\partial}
\newcommand\cl{\overline}
\newcommand\tsb[1]{_{\text{#1}}}
\newcommand\set[1]{\{#1\}}
\newcommand\defset[3][\colon]{\set{#2#1#3}}
\newcommand\R{\mathbb R}
\renewcommand\P{\mathcal P}
\newcommand\T{\mathcal T}
\newcommand\U{\mathcal U}
\newcommand\V{\mathcal V}
\newcommand\X{\mathcal X}
\newcommand\wmax{w\tsb{max}}
\let\0\emptyset
\DeclareMathOperator\tr{tr}

\begin{document}

\title{Modeling and numerical simulation of
fully Eulerian fluid-structure interaction
using cut finite elements}


\author{Stefan Frei~\thanks{
Stefan Frei, Department of Mathematics \& Statistics, University of Konstanz, {stefan.frei@uni-konstanz.de}}\, , 
Tobias Knoke, 
Marc C. Steinbach,\\ 
Anne-Kathrin Wenske, 
Thomas Wick~\thanks{Tobias Knoke, Marc C. Steinbach, Anne-Kathrin Wenske,
Thomas Wick, Leibniz University Hannover, Institute of Applied Mathematics {\{knoke,mcs,wenske,thomas.wick\}@ifam.uni-hannover.de}}
}
%
%
\maketitle

\begin{abstract}
We present a monolithic finite element formulation
  for (nonlinear) fluid-structure interaction in Eulerian coordinates.
  For the discretization we employ an unfitted finite element method
  based on inf-sup stable finite elements.
  So-called ghost penalty terms are used
  to guarantee the robustness of the approach
  independently of the way the interface cuts the finite element mesh.
  The resulting system is solved in a monolithic fashion
  using Newton's method.
  Our developments are tested on a numerical example
  with fixed interface.
\end{abstract}

\section{Introduction}
\label{sec:1}
In this work, we investigate a cut finite element discretization
for fully Eulerian fluid-structure interaction (FSI).
In contrast to an \textit{Arbitrary Lagrangian Eulerian} (ALE)
approach \cite{DoneaSurvey, Hughes1981},
the benefit of a fully Eulerian formulation for FSI
lies in its ability to handle (very) large deformations,
topology changes, and contact problems in a straight-forward way,
see e.g.\
\cite{FrRiWi16_JCP,FrRi17,HechtPironneau2016,BurmanFernandezFrei2020}.

The fully Eulerian approach for fluid-structure interaction
has been introduced in~\cite{Du06,CoMaMi08} and has since then been investigated and improved in several studies, such as \cite{richterWick2010,SuIiTaTaMa11,Wi12_fsi_euler,RICHTER2013227,LaRuiQua13,
Sun20141,RATH2023112188}. The idea is to formulate both the flow
and the solid problem in Eulerian coordinates in time-dependent domains $\Omega_f(t)$ resp.~$\Omega_s(t)$. An accurate numerical method requires the resolution of the interface $\Gamma(t)$ separating $\Omega_f(t)$ and $\Omega_s(t)$, which can move freely depending on the solid displacements.
The construction of a fitted finite element method is cumbersome when the interface moves, see e.g.~\cite{FreiRichter2014, FreiPressure2019}. An elegant alternative is given by the cut finite element method~\cite{HANSBO2002elliptic, hansbo2005Nitsche, Burman2015CutFEM} which is based on a fixed finite
element mesh for all times. To our knowledge, this approach has not been used before in the context of fully Eulerian fluid-structure interaction. The present work is thus a first step towards such an unfitted fully Eulerian FSI formulation. As a starting point we concentrate on fixed interfaces
(meaning infinitesimal displacements) in this work.

In cut finite element methods (CutFEM) \cite{Burman2015CutFEM}, interface conditions are imposed by means of Nitsche's method \cite{Nitsche1971},
see also \cite{HANSBO2002elliptic, hansbo2005Nitsche}. Moreover, additional stabilization via ghost penalty terms at the faces of cut cells is proposed, since the condition number of the system matrix suffers from cells cut into vastly
different sizes, see also \cite{BURMAN2010GhostPenalty, BURMAN2012FictitiousDomainII}. This adaptation of Nitsche's method is used in \cite{BURMAN2014FSI}, where
linear Stokes flow is coupled to a linear elastic structure through separate overlapping meshes, where the solid is described in Lagrangian coordinates on a fitted mesh and glued to the (unfitted) fluid mesh. Fictitious domain methods using cut elements with stabilized Nitsche's method for Stokes’ problem
were investigated in \cite{burman_hansbo_2014}.

In this work, we employ a cut finite element method for realizing a variational-monolithic fully Eulerian fluid-structure interaction formulation
on a fixed single mesh. This
is specifically in extension to \cite{BURMAN2014FSI} and \cite{burman_hansbo_2014}
from the ghost penalty viewpoint and all previously mentioned fully Eulerian fluid-structure
interaction references.
We propose a new weight function to balance
the ghost penalty terms with respect to the cuts. Taylor-Hood elements are employed for
the spatial discretization and a backward Euler scheme for temporal discretization.
The resulting discrete monolithic formulation
is treated all-at-once in the linear and nonlinear solvers.

The outline of this paper is as follows. In \cref{sec_FSI} our
fully Eulerian fluid-structure interaction formulation is presented.
Then, in \cref{sec_nonlinear}, the cut finite element discretization
and ghost penalty terms are introduced and described in detail.
Finally, in \cref{sec_tests}, numerical simulations are carried out,
including numerical convergence studies and a comparison
to computations using the ALE method.

\section{Fluid-structure interaction system}
\label{sec_FSI}

\subsection{Strong form}
Let $\Omega \subset \R^d$ be a bounded domain
with $d=2$, which is partitioned into a
(fixed) fluid subdomain $\Omega_f$ and a
(fixed) solid subdomain $\Omega_s$ such that
$\cl\Omega = \cl\Omega_f \cup \cl\Omega_s$ with
$\Omega_f \cap \Omega_s = \0$.
We assume that both $\Omega_f$ and $\Omega_s$
are parameterized by a $C^{1,1}$ boundary,
such that all terms arising in the following equations are well-defined.
Next, let $\Gamma_i := \cl\Omega_f \cap \cl\Omega_s$
be the interface between the subdomains, and
$\Gamma_f^D \subset \Gamma_f := \bd\Omega_f \cap \bd\Omega$ and
$\Gamma_s^D \subset \Gamma_s := \bd\Omega_s \cap \bd\Omega$.
For the fluid velocity $v_f\colon \cl\Omega_f \times [0, T] \to \R^d$,
the pressure $p\colon \cl\Omega_f \times [0, T] \to \R$,
the solid velocity $v_s\colon \cl\Omega_s \times [0, T] \to \R^d$ and
the displacement $u\colon \cl\Omega_s \times [0, T] \to \R^d$,
we define stresses
$\sigma_f := \sigma_f(v_f,p) := \rho_f \nu_f (\nabla v_f + \nabla v_f^\top) - pI$
and $\sigma_s := 2 \mu_s E_s + \lambda_s \tr(E_s)I$, where
$E_s := \frac12 (\nabla u + \nabla u^\top + \nabla u^\top \cdot \nabla u)$
denotes the nonlinear Green-Lagrange strain,
$\mu_s$ and $\lambda_s$ are the Lam\'{e} parameters,
$\rho_f$ and $\rho_s$ are the densities of the fluid and the solid,
and $\nu_f$ is the fluid viscosity.
Moreover, $f\colon \Omega \times [0, T] \to \R^d$
is a given right-hand side function,
$v_f^D\colon \Gamma_f^D \times [0, T] \to \R^d$
and
$u^D\colon \Gamma_s^D \times [0, T] \to \R^d$
are functions on the Dirichlet boundaries, and
$\smash{v_f^0}\colon \Omega_f \to \R^d$,
$v_s^0\colon \Omega_s \to \R^d$ and
$u^0\colon \Omega_s \to \R^d$ finally describe initial values.
The fully Eulerian FSI system is then defined as follows:
Find $(v_f,p,v_s,u)$ such that
\begin{align*}
&\left\{\begin{array}{lll}
        \rho_f \partial_t v_f+\rho_f (v_f \cdot \nabla)v_f-\nabla \cdot \sigma_f &=\rho_f f \quad & \text{in } \Omega_f \times (0, T), \\
        \nabla \cdot v_f &=0 \quad & \text{in } \Omega_f \times (0, T), \\
        v_f &=v_f^D \quad & \text{on } \Gamma_f^D \times (0, T), \\
        \rho_f \nu_f \partial_n v_f-p n &= 0
        & \text{on } \Gamma_f \setminus \Gamma_f^D \times (0, T), \\
        v_f &= v_f^0 & \text{in } \Omega_f \times \set{0},
       \end{array} \right. \\[-0.5pt]
&\left\{\begin{array}{lll}
        \rho_s \partial_t v_s +\rho_s (v_s \cdot \nabla)v_s-\nabla \cdot \sigma_s &=\rho_s f \quad & \text{in } \Omega_s \times (0, T), \\
        \partial_t u + (v_s \cdot \nabla)u - v_s &=0 \quad & \text{in } \Omega_s \times (0, T), \\
        u &=u^D \quad & \text{on } \Gamma_s^D \times (0, T), \\
        \sigma_s \cdot n &=0 \quad & \text{on } \Gamma_s \setminus \Gamma_s^D \times (0, T), \\
        u &= u^0 & \text{in } \Omega_s \times \set{0}, \\
        v_s &= v_s^0 & \text{in } \Omega_s \times \set{0},
       \end{array} \right. \\[-0.5pt]
&\left\{\begin{array}{lll}
        v_f &= v_s \quad & \text{on } \Gamma_i \times (0, T), \\
        \sigma_f \cdot n &=\sigma_s \cdot n \quad & \text{on }
\Gamma_i \times (0, T).
       \end{array} \right.
\end{align*}

\subsection{Weak formulation}

Let $\V_f := H_0^1(\Omega_f;\Gamma_f^D)$,
$\V_s := H^1(\Omega_s)$, 
$\U := H_0^1(\Omega_s;\Gamma_s^D)$ and
$\P := L^2(\Omega_f)$ be given function spaces.
Here, $\P$ is sufficient for a unique pressure due to the outflow condition on
$\Gamma_f \setminus \Gamma_f^D$.
The product space is defined as
$\X := \V_f \times \V_s \times \U \times \P$.
\begin{Problem}
Find $v_f \in v_f^D + \V_f$, $p \in \P$, $v_s \in \V_s$
and $u \in u^D + \U$ such that
$v_f = v_s$ on $\Gamma_i$ and for all $(\phi_f, \psi, \phi_s, \xi) \in \X$:
\begin{align*}
  &\rho_f \prodf{\partial_t v_f}{\phi_f} +
    \rho_f \prodf{v_f \cdot \nabla v_f}{\phi_f} +
    \prodf{\sigma_f}{\nabla \phi_f} +
    \prodf{\nabla \cdot v_f}{\xi} \\[-2.5pt]
  &- \prodd{\rho_f \nu_f \nabla v_f^\top n_f}{\phi_f} +
    \rho_s \prods{\partial_t v_s + v_s \cdot \nabla v_s}{\phi_s} +
    \prods{\sigma_s}{\nabla \phi_s} \\[-3pt]
  &+\prods{\partial_t u + v_s \cdot \nabla u - v_s}{\psi}
    - \prodi{\sigma_f \cdot n_s}{\phi_f - \phi_s}
    =
    \rho_f \prodf{f}{\phi_f} + \rho_s \prods{f}{\phi_s}
    .
\end{align*}
\end{Problem}

\subsection{Discretization and ghost penalties}

To discretize in time we apply the backward Euler method.
For spatial discretization we use
continuous quadratic elements for the fluid velocity and
continuous linear elements for the remaining solution components.
Let $\T_h$ be a quasi-uniform triangulation of $\Omega$
that is fitted to the boundary of the domain $\Omega$
but not to the interface $\Gamma_i$,
where $\Gamma_i$ is described by a level set function.
Moreover, let
\begin{align*}
  \T_h^f := \defset{T \in \T_h}{T \cap \Omega_f \ne \0}
  \quad\text{and}\quad
  \T_h^s := \defset{T \in \T_h}{T \cap \Omega_s \ne \0}
\end{align*}
be overlapping sub-triangulations.
We use the following finite element spaces on $\T_h^i$:
\begin{align*}
  V_{h,i}^{(r)} :=
  \defset{\phi \in C(\cl{\Omega_h^i})}
  {\phi|_T \in P_r(T) \ \forall T \in \T_h^i},
  \quad i \in \set{f, s},
\end{align*}
where $\Omega_h^i$ denotes the domain spanned by the cells $T\in {\cal T}_h^i$.
We define $\V_{h,f}:=V_{h,f}^{(2)} \cap H_0^1(\Omega_f;\Gamma_f^D)$,
$\V_{h,s} :=V_{h,s}^{(1)}$,
$\U_h := V_{h,s}^{(1)}\cap H_0^1(\Omega_s;\Gamma_s^D)$,
$\P_h := V_{h,f}^{(1)}$ and
$\X_h := \V_{h,f} \times \V_{h,s} \times \U_h \times \P_h$.

The interface conditions are then imposed by additional terms:
the Nitsche terms \eqref{pb2_4} which ensure
that the interface condition $v_f = v_s$ is satisfied,
a stabilization term \eqref{pb2_5}
to control the pressure (see \cite{BURMAN2014FSI}),
and ghost penalty terms \eqref{pb2_6},\eqref{pb2_8},
described by the ghost penalty functions $g^{h,w}(\,\cdot\,,\,\cdot\,)$
around the interface zone, that
extend the coercivity of the bilinear form over the interface cells and
increase stability.
\begin{Problem}
The discrete weak formulation reads: For $n=1,\ldots,N$ find
fluid velocity, pressure, solid velocity and displacement
$(v_f^h, p^h, v_s^h, u^h) := (v_f^{h, n}, p^{h, n}, v_s^{h, n}, u^{h, n})
\in \{v_f^D,0,u^D,0\} + \X_h$, where
$(v_f^{h, n-1}, p^{h, n-1}, v_s^{h, n-1}, u^{h, n-1})$
are the solutions of the previous time step, such that
for all $(\phi_f^h, \psi^h, \phi_s^h, \xi^h) \in \X_h$:
\begin{align}
  &\rho_f \prodf{v_f^h}{\phi_f^h} +
    \rho_f k\prodf{ v_f^h \cdot \nabla v_f^h}{\phi_f^h} +
    k \prodf{\sigma_f^h}{\nabla \phi_f^h} +
    \prodf{\nabla \cdot v_f^h}{\xi^h}  \label{pb2_1} \\[-.5pt]
  &\quad- k\prodd{ \rho_f \nu_f (\nabla v_f^h)^\top n_f}{\phi_f^h}
    + \rho_s \prods{v_s^h}{\phi_s^h}
    + k \prods{\sigma_s^h}{\nabla \phi_s^h} \label{pb2_2} \\[-.5pt]
  &\quad+ \rho_s k\prods{ v_s^h \cdot \nabla v_s^h}{\phi_s^h}
    + \prods{u^h +k(v_s^h \cdot \nabla u^h - v_s^h)}{\psi^h}
    \label{pb2_3} \\[-.5pt]
  &\quad+ k h^{-1} \rho_f \nu_f \gamma_N
    \prodi{v_f^h - v_s^h}{\phi_f^h - \phi_s^h}
    - k \prodi{\sigma_f^h \cdot n_f}{\phi_f^h - \phi_s^h}
    \label{pb2_4} \\[-.5pt]
  &\quad- k \prodi{v_f^h - v_s^h}{\sigma_f^h (\phi_f^h, -\xi^h) \cdot n_f}
    \label{pb2_5} \\[-.5pt]
  &\quad+ 2 \rho_f \nu_f k g_{v_f}^h(v_f^h, \phi_f^h) +
    \rho_s g_{v_s}^h(v_s^h, \phi_s^h) +
    k g_p^h(p^h, \xi^h) +
    2 \mu_s k g_u^h(u^h, \phi_s^h) \label{pb2_6} \\[-.5pt]
  &= \rho_f k \prodf{f}{\phi_f^h} +
    \rho_f \prodf{v_f^{h, n-1}}{\phi_f^h} +
    \rho_s k \prods{ f}{\phi_s^h} +
    \rho_s \prods{v_s^{h, n-1}}{\phi_s^h} \label{pb2_7} \\[-.5pt]
  &\quad+ \prods{u^{h, n-1}}{\psi^h}
    + g_{v_s}^{h, w}(v_s^{h, n-1}, \phi_s^h) \label{pb2_8}
      .
\end{align}
Here $k=t_n - t_{n-1}>0$ is the time step size, $h>0$ the spatial discretization
parameter, namely the maximum element size,
and $\gamma_N>0$ denotes the Nitsche parameter.
\end{Problem}
As usual, we express the weak form more compactly
in terms of a semi-linear form:
Find $U_h^n := (v_f^{h, n},v_s^{h, n},u^{h, n},p^{h, n})
\in \{v_f^D,0,u^D,0\} + \X_h$
for the time steps $n = 1, \dots, N$,
such that
$
  A(U_h^n)(\Psi) = F(\Psi)
$
for all $\Psi \in \X_h$.

Let $\mathcal{F}_G^f$ denote the set of element faces
$F = \bar{K}_1 \cap \bar{K}_2$ of the triangulation $\T_h^f$
that do not lie on the boundary $\bd\Omega$,
such that at least one of the cells $K_j$ is intersected by the interface
$(K_j \cap \Gamma_i \ne \0$ for $j \in \set{1, 2}$).
Analogously, we define $\mathcal{F}_G^s$
as the set of corresponding faces of the triangulation $\T_h^s$.
For a cell $K$ cut by the interface we denote by $K\tsb{in}$
the part of the cell inside the considered subdomain.
Using the jump terms
$\jump{u} = u_1|_\Gamma - u_2|_\Gamma$ with $u_i = u|_{\Omega_i}$
and $\Jump{i} := \iprod[F]{\jumpn1i}{\jumpn2i}$,
the ghost penalty functions (with parameters
$\gamma_{v_f}$, $\gamma_p$, $\gamma_{v_s}$, $\gamma_u$)
are defined as follows:
\newcommand\sumFf{\sum_{F \in \mathcal{F}_G^f}}
\newcommand\sumFKf{\sumFf \sum_{K\mathpunct: F \in \cl K}}
\newcommand\sumFs{\sum_{F \in \mathcal{F}_G^s}}
\newcommand\sumFKs{\sumFs \sum_{K\mathpunct: F \in \cl K}}
\begin{align*}
  g_{v_f}^{h, w}(\varphi_1, \varphi_2)
  &:= \gamma_{v_f} \sumFKf
    w(\kappa_K) \Bigl( h \Jump1 + \frac{h^3}{4} \Jump2 \Bigr), \\
  g_p^{h, w}(\varphi_1, \varphi_2)
  &:= \gamma_p \sumFKf
  w(\kappa_K) h^3 \Jump1, \\
  g_{v_s}^{h, w}(\varphi_1, \varphi_2)
  &:= \gamma_{v_s} \sumFKs w(\kappa_K) h^3 \Jump1, \\
  g_u^{h, w}(\varphi_1, \varphi_2)
  &:= \gamma_u \sumFKs w(\kappa_K) h \Jump1.
\end{align*}
Here we apply a novel weight function
$w\colon [0,1] \to [\frac12 \wmax^{-1},\frac12 \wmax]$,
$\kappa \mapsto \frac12 \wmax^{1 - 2 \kappa}$, with $\wmax \ge 1$,
which scales the ghost penalties dependent on the cell cuts by taking
the portion of the inside cell part,
$\kappa_K := \text{meas}(K_\text{in}) / \text{meas}(K)$, as the argument.
Thus we penalize ``bad cuts'' more severely
while ``good cuts'' (where a sufficiently large portion of the cell lies inside)
are penalized less severely.
Moreover,
  the conventional ghost penalty terms are recovered as the special case
  where $\wmax = 1$, hence $w \equiv 0.5$.

\section{Nonlinear solution}
\label{sec_nonlinear}
To employ Newton's method, we need the derivative
\begin{align*}
  A'&(U_h^{n,j})(\delta U_h, \Psi) =
      \lim_{\epsilon \to 0} \frac{1}{\epsilon} \left(
      A(U_h^{n,j} + \epsilon \delta U_h)(\Psi) - A(U_h^{n,j})(\Psi) \right) \\
    &= \rho_f \prodf{\delta v_f^h}{\phi_f^h}
      + \rho_f k \prodf{\delta v_f^h \cdot \nabla v_f^{h,j} +
      v_f^{h,j} \cdot \nabla \delta v_f^h}{\phi_f^h}\\
    &+ k \prodf{\sigma_f^h (\delta v_f^h, \delta p^h)}{\nabla \phi_f^h}
      + \prodf{\nabla \cdot \delta v_f^h}{\xi^h}
      - k \rho_f \nu_f
      \smash{\prodd{\nabla {\delta v_f^h}^\top n_f}{\phi_f^h}} \\
    &+ \rho_s \prods{\delta v_s^h}{\phi_s^h}
      + \rho_s k \prods{\delta v_s^h \cdot \nabla v_s^{h,j} +
      v_s^{h,j} \cdot \nabla \delta v_s^h}{\phi_s^h} \\
    &+ k \prods{{\sigma_s^h}'(u^{h,j})(\delta u^h)}{\nabla \phi_s^h}
      + 
      \prods{\delta u^h + k(\delta v_s^h \cdot \nabla u^{h,j} +
      v_s^{h,j} \cdot \nabla \delta u^h - \delta v_s^h)}{\psi^h} \\
    &+ k h^{-1} \rho_f \nu_f \gamma_N
      \prodi{\delta v_f^h - \delta v_s^h}{\phi_f^h - \phi_s^h}
      - k \prodi{\sigma_f^h (\delta v_f^h, \delta p^h) \cdot n_f}
      {\phi_f^h - \phi_s^h} \\
    &- k \prodi{\delta v_f^h - \delta v_s^h}
      {\sigma_f^h (\phi_f^h, -\xi^h) \cdot n_f} \\
    &+ 2 \rho_f \nu_f k g_{v_f}^h (\delta v_f^h, \phi_f^h)
      + \rho_s g_{v_s}^h (\delta v_s^h, \phi_s^h)
      + k g_p^h (\delta p^h, \xi^h)
      + 2 \mu_s k g_u^h (\delta u^h, \phi_s^h).
\end{align*}
where
\begin{align*}
  {\sigma_s^h}'(u^{h,j})(\delta u^h)
  &:= 2 \mu_s E_s'(u^{h,j})(\delta u^h) +
    \lambda_s \tr \bigl(E_s'(u^{h,j})(\delta u^h)\bigr), \\[-1pt]
  E_s'(u^{h,j})(\delta u^h)
  &:= \tfrac12 \bigl( \nabla \delta u^h + (\nabla \delta u^h)^\top
    + (\nabla \delta u^h)^\top \cdot \nabla u^{j,h}
    + (\nabla u^{j,h})^\top \cdot \nabla \delta u^h \bigr).
\end{align*}

With the step length $\alpha^j \in (0,1]$ determined by a line search,
Newton's method then takes the following form:
Given an initial guess $U_h^{n,0} \in \set{\smash{v_f^D},0,u^D,0} + \X_h$,
such as $U_h^{n,0} := U_n^{n-1}$, find $\delta U_h \in \X_h$
for $j = 0, 1, 2, \dots$, such that for all $\Psi \in \X_h$:
\begin{align*}
  A'(U_h^{n,j})(\delta U_h,\Psi)
  &= -A(U_h^{n,j})(\Psi) + F(\Psi), \quad
  U_h^{n,j+1} = U_h^{n,j} + \alpha^j \delta U_h.
\end{align*}

\section{Numerical test: modified ``flow around a cylinder benchmark''}
\label{sec_tests}
In this section, we apply our numerical framework to a model problem
inspired by the \textit{flow around a cylinder benchmark}~\cite{schafer1996}.
We use Newton's method for the nonlinear solution, and therein for the
linear systems the parallel sparse solver MUMPS \cite{Amestoy2001MUMPS}.
The implementation is based on the open-source finite element
library deal.II \cite{dealII94}, in particular step 85 of the tutorial programs.
Comparative computations with an arbitrary Eulerian-Lagrangian
fluid-structure interaction formulation are also performed with the open-source
code \cite{Wi13_fsi_with_deal}.
For our computations we neglect the convection terms
$\rho_s (v_s \cdot \nabla)v_s$ and $(v_s \cdot \nabla)u$ in the structure.

\begin{figure}[tp]
\centering
\resizebox{\textwidth}{!}{%
	\begin{tikzpicture}[scale = 4.5]
		\draw (0,0) to node [below, fill = none] {$\Gamma\tsb{wall}$} (2.2, 0)
			to node [right, fill = none] {$\Gamma\tsb{out}$} (2.2, 0.41)
			to node [above, fill = none] {$\Gamma\tsb{wall}$} (0, 0.41)
			to node [left, fill = none] {$\Gamma\tsb{in}$} (0,0);
		\filldraw[pattern = north east lines] (0.2, 0.2) circle (0.05);
		\draw[fill=white] (0.2, 0.2) circle (0.01);
		\draw[-latex] (0, 0.12) -- (0.065, 0.12);
		\draw[-latex] (0, 0.205) -- (0.115, 0.205);
		\draw[-latex] (0, 0.29) -- (0.065, 0.29);
		\node (inner_wall) at (0.5, 0.33) {$\Gamma\tsb{wall}$};
		\draw[-latex] (inner_wall) -- (0.209, 0.2044);
		\node[below left] at (0,0) {(0,0)};
		\node[above left] at (0,0.41) {(0,0.41)};
		\node[above right] at (2.2, 0.41) {(2.2,0.41)};
		\node[below right] at (2.2, 0) {(2.2,0)};
	\end{tikzpicture}
	}
\caption{Configuration of laminar flow around an elastic shell with center at
  (0.2, 0.2), inner radius of 0.01 and outer radius of 0.05.}
\label{fig:channel}
\end{figure}
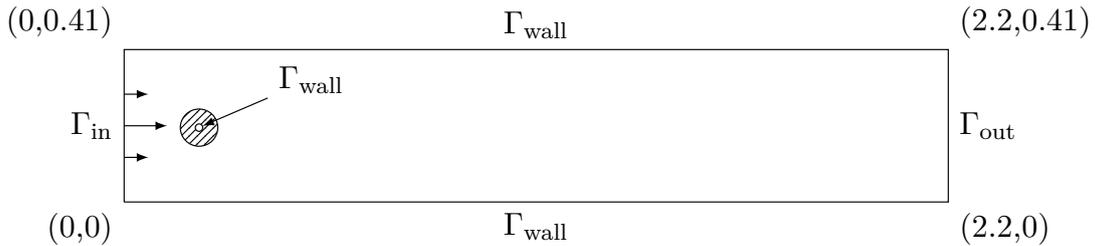

In our modification of the laminar flow benchmark
\cite{schafer1996} the cylindrical hole
is replaced by an elastic solid with a hole in the middle
as depicted in \cref{fig:channel}.
The remaining channel is filled with an incompressible Newtonian fluid.
At the boundary $\Gamma\tsb{in}$, we impose a parabolic inflow profile given by
\begin{equation*}
  v_f(0,y) = 1.5 \bar{U} (4y(0.41-y))/0.41^2
\end{equation*}
with mean velocity $\bar{U} = \SI{0.2}{ms^{-1}}$.
At $\Gamma\tsb{out}$ a do-nothing outflow condition is applied.
The boundaries $\Gamma\tsb{wall}$ supply
a no-slip-condition for the fluid and
a homogenous Dirichlet condition for the solid deformation.
We start the time-stepping with homogeneous initial conditions
and increase the inflow gradually by setting
\begin{equation*}
  v_f(t,0,y) = \tfrac12 (1-\cos(\tfrac{t\pi}{2})) v_f(0,y)
  \quad\text{for } t < 2.
\end{equation*}

The material parameters are based on the FSI-1 benchmark \cite{Turek2006FSIBM}:
fluid density $\rho_f = \SI{1000}{kg/m^3}$,
fluid viscosity $\nu_f = \SI{0.001}{m^2/s}$,
solid density $\rho_s = \SI{1000}{kg/m^3}$,
Lamé coefficients \num{0.5e+6} and \num{2.0e+6}.

The mesh consists of rectangular elements with shape regular cells
except for the area of the circular solid domain $\Omega_s$
and a small neighborhood thereof.

\begin{figure}[t]
  \centering
  \def\plot#1{\includegraphics[width=0.9\textwidth,
    trim=1.2cm 20cm 6cm 6cm,clip]{eulerian_#1.png}}
  \plot{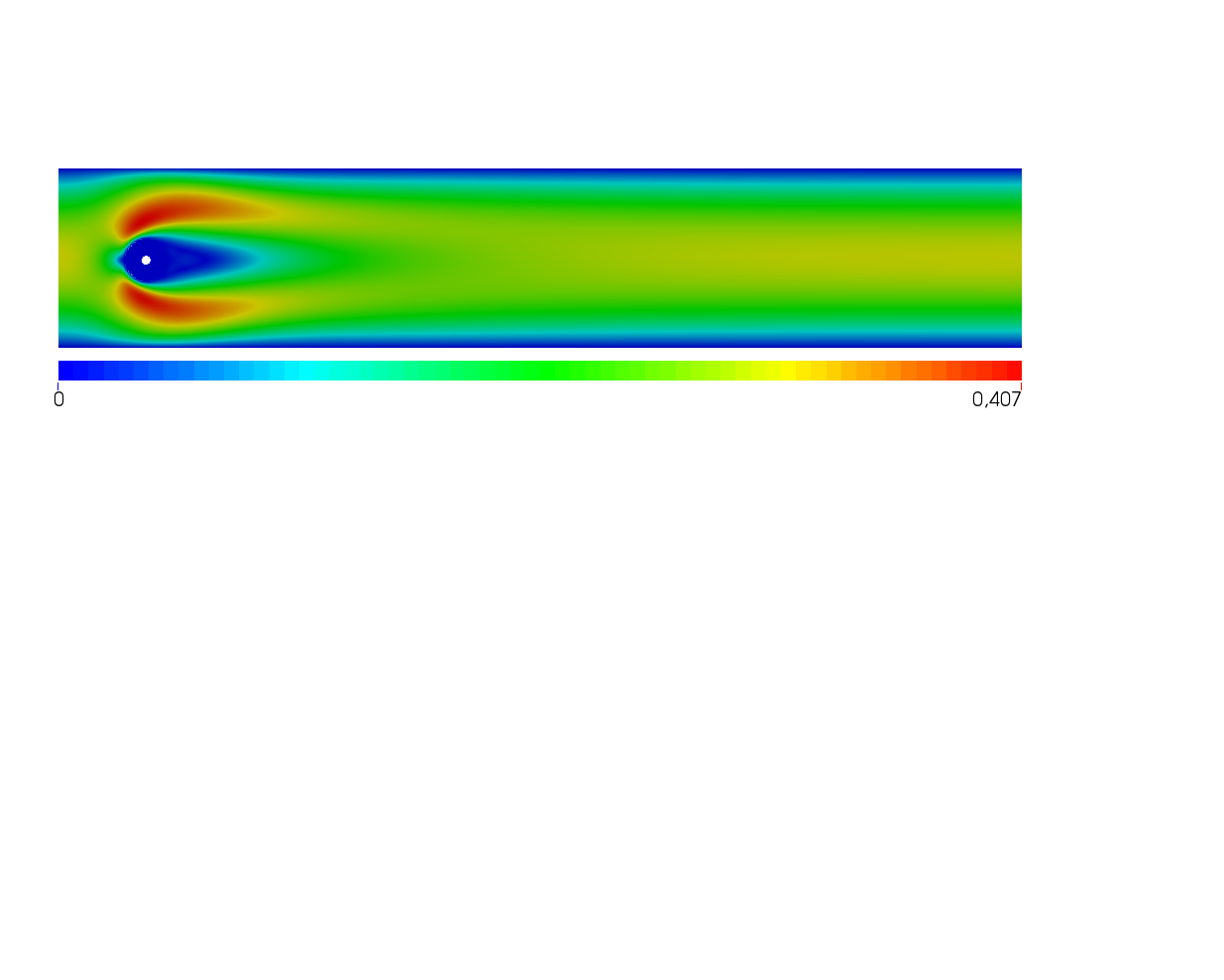}
  \plot{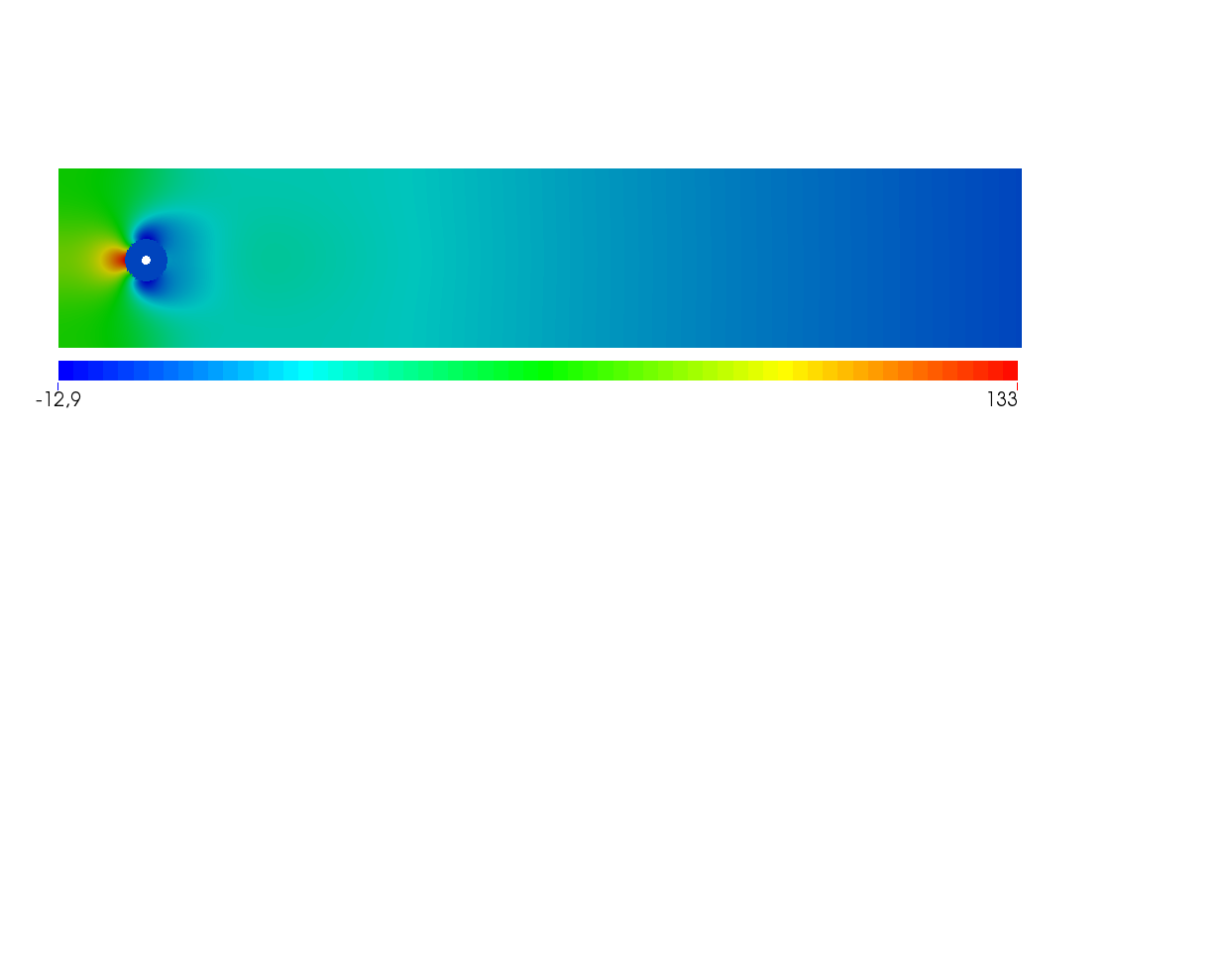}
  \plot{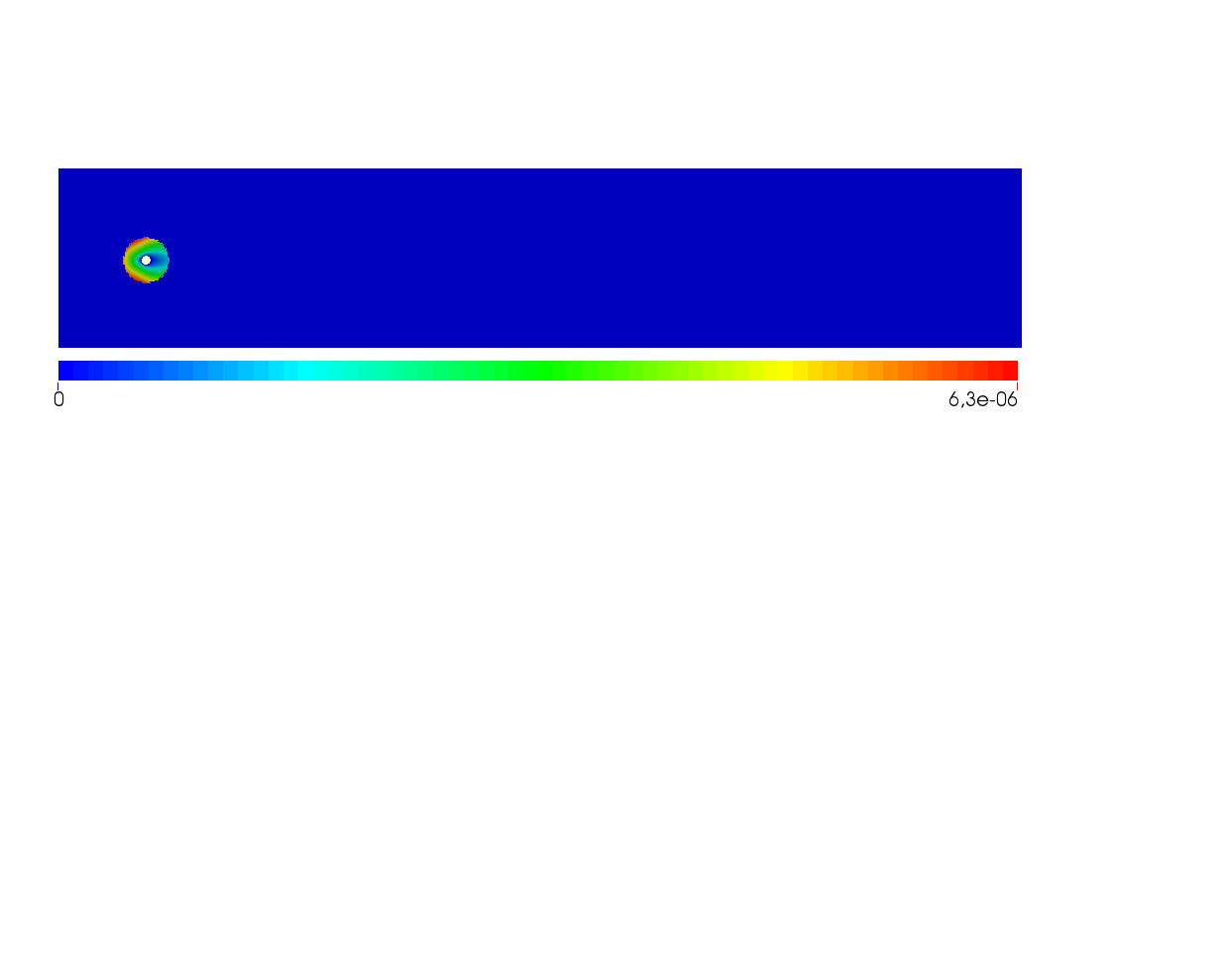}
  \caption{From top to bottom: Fluid speed $\norm{v_f}_2$,
    pressure and displacement using the fully Eulerian approach
    at time $T = 25$ with $\wmax = 3$.
    As the results using the ALE approach are visually very similar
    to the fully Eulerian approach, we refrain from including their plots.}
  \label{fig:euler-1}
\end{figure}

We simulate over the time interval $[0, T]$ with $T = 25$
using the time step size $\Delta t = 1.0$.
We choose the Nitsche parameter $\gamma_N = 10$,
the two weight parameters $\wmax = 3$ and $\wmax = 1$,
and ghost penalty parameters
$\gamma_{v_f} = \gamma_p = \gamma_{v_s} = \gamma_u = 10^{-3}$.
For Newton's method we use the absolute tolerance $10^{-8}$.

The interface between fluid and solid
is artificially fixed to enable simple implicit time stepping.
The solid material is comparatively stiff,
resulting in small deformations.
As quantities of interest we choose the fluid velocity
at the center point of the outflow boundary, $(2.2, 0.205)$,
as well as drag and lift forces around the solid.
The latter are given by the line integral over the interface $\Gamma_i$,
$
  (F_D, F_L) = \int_{\Gamma_i} \sigma_f \cdot n_f,
$
where $n_f$ is the normal vector on the interface
pointing towards the fluid domain.
The results are shown in \cref{tab:euler-1,fig:euler-1}.

\begin{table}[t]
  \def\0{\hphantom0}
  \centering
  \caption{Quantities of interest using the FSI-1 parameters
    on uniformly refined meshes.
    The results in the top half correspond to the weight $\wmax=3$
    and at the bottom to $\wmax=1$.
    Here, $L$ denotes the refinement level.
    The values in the $L^2$ norm are taken at time $T = 25$.
    The drag and lift forces as well as the fluid velocity at the outflow
    are averaged over the time interval.
  }
  \label{tab:euler-1}
  \begin{tabular}{c@{\quad}r@{\quad}c@{\quad}c@{\quad}c@{\quad}c@{\quad}c@{\quad}c}
    \toprule
    $L$  & \# dofs   & $\normL{\nabla v_f}$ & $\normL{p}$ &  $\normL{\nabla u}$ & $F_D$   & $F_L$  & $v_f(2.2,0.205)$\\
    & & & $[\times 10^{-8}]$ & $[\times 10^{-5}]$ & & & $([\times 10^{-1}],[\times 10^{-4}])$\\
    \midrule
    0 & 3740   & 2.4152                 & 24.4730                       & 1.127               &\08.6278  & 0.0323  & $(2.734295, -1.5671)$\\
    1 & 13944  & 2.4251                 & 24.7050                       & 1.059               & 10.3264 & 0.0193  & $(2.734141, -1.6096)$\\
    2 & 53972  & 2.4283                 & 24.7473                       & 1.067               & 10.5273 & 0.0209  & $(2.734137, -1.6148)$\\
    3 & 211980 & 2.4294                 & 24.7572                       & 1.081               & 10.5484 & 0.0215  & $(2.734135, -1.6166)$\\
    4 & 839900 & 2.4297                 & 24.7587                       & 1.088               & 10.5654 & 0.0215  & $(2.734136, -1.6170)$\\
    \midrule
    0 & 3740   & 2.4153                 & 24.4752                       & 1.127               &\08.6236  & 0.0322 & $(2.734302, -1.5681)$\\
    1 & 13944  & 2.4251                 & 24.7058                       & 1.060               & 10.3357 & 0.0192 & $(2.734144, -1.6097)$\\
    2 & 53972  & 2.4284                 & 24.7487                       & 1.067               & 10.5495 & 0.0212 & $(2.734133, -1.6150)$\\
    3 & 211980 & 2.4295                 & 24.7577                       & 1.081               & 10.5709 & 0.0215 & $(2.734135, -1.6167)$\\
    4 & 839900 & 2.4297                 & 24.7589                       & 1.088               & 10.5864 & 0.0216 & $(2.734136, -1.6170)$\\
    \bottomrule
  \end{tabular}
\end{table}

We observe no qualitative difference between our choices of the weights.
However, we note that the fluctuation of the fluid velocity
in cut cells is less severe for $\wmax = 3$
as compared to the traditional ghost penalization.

We compare these results with a corresponding computation using an
ALE approach~\cite{Wi13_fsi_with_deal}.
Here, the mesh is fitted to the interface
and quadratic elements are used for the structural variables,
which implies that the number of
degrees of freedom differs from the Eulerian approach.
The ALE results are shown in \cref{tab:ale-1}.

\begin{table}[t]
  \caption{Quantities of interest using the ALE approach
    as with \cref{tab:euler-1}.}
  \label{tab:ale-1}
  \def\0{\hphantom0}
  \centering
  \begin{tabular}{c@{\quad}r@{\quad}c@{\quad}c@{\quad}c@{\quad}c@{\quad}c@{\quad}c}
    \toprule
    $L$   & \# dofs   & $\normL{\nabla v_f}$ & $\normL{p}$  & $\normL{\nabla u}$ & $F_D$   & $F_L$  & $v_f(2.2,0.205)$\\
    & & & $[\times 10^{-8}]$ & $[\times 10^{-5}]$ & & & $([\times 10^{-1}],[\times 10^{-4}])$\\
    \midrule
    0 & 8080   & 2.4005                 & 24.4431                       & 1.064               &\09.6868  & 0.1043 & $(2.73447, -1.56003)$\\
    1 & 31360  & 2.4183                 & 24.6757                       & 1.078               & 10.3359 & 0.0323 & $(2.73421, -1.59833)$\\
    2 & 123520 & 2.4210                 & 24.7056                       & 1.081               & 10.4933 & 0.0220 & $(2.73418, -1.60355)$\\
    3 & 490240 & 2.4211                 & 24.7058                       & 1.082               & 10.5150 & 0.0205 & $(2.73418, -1.60389)$\\
    \bottomrule
  \end{tabular}
\end{table}

The two approaches are in good agreement, as the above tables show.
We observe convergence for $h \to 0$ in all quantities of interest.
The small deviations between ALE and fully Eulerian computations
can be explained by the time discretization errors,
as the time step $k$ is fixed.
To further compare the two solutions we investigate the values
of fluid pressure and its speed along three vertical lines
$\defset{(x,y)}{0 \le y \le 0.41}$ for $x \in \set{0.15,0.25,2.2}$.
\Cref{fig:plots1} depicts these results at the final time $T = 25$.
Again, we observe a generally good agreement
between our solution and the ALE model.

\pgfplotsset{compat=1.14,
  width=180pt,height=141pt,grid style={gray!15},grid=both,
  xmin=-0.03,xmax=0.44,minor x tick num=1,minor y tick num=1,
}
\begin{figure}[t]
\centering
  \def\plot#1#2{%
    \begin{tikzpicture}
      \footnotesize
      \begin{axis}
        \draw[gray!40] (-1,0)--(1,0);
        \addplot[blue,only marks,mark size=1,mark=x]
        table[y index=#2] {line_#1_ale};
        \addplot[red,only marks,mark size=1,mark=o]
        table[y index=#2] {line_#1_eulerian};
      \end{axis}
    \end{tikzpicture}
  }
  \plot12\hfil\plot11\\
  \plot22\hfil\plot21\\
  \plot32\hfil\plot31
  \caption{Profiles of fluid speed $\norm{v_f}_2$ (left) and pressure (right)
    at time $T = 25$ with $\wmax=3$ along three vertical lines:
    at $x = 0.15$ in front of the solid (top),
    at $x = 0.25$ behind the solid (middle),
    and at $x = 2.2$ at the outflow boundary (bottom);
    ALE = $\color{blue}\times$,
    Euler = \raisebox{-.2ex}{\large\color{red}$\circ$}.
    }
  \label{fig:plots1}
\end{figure}
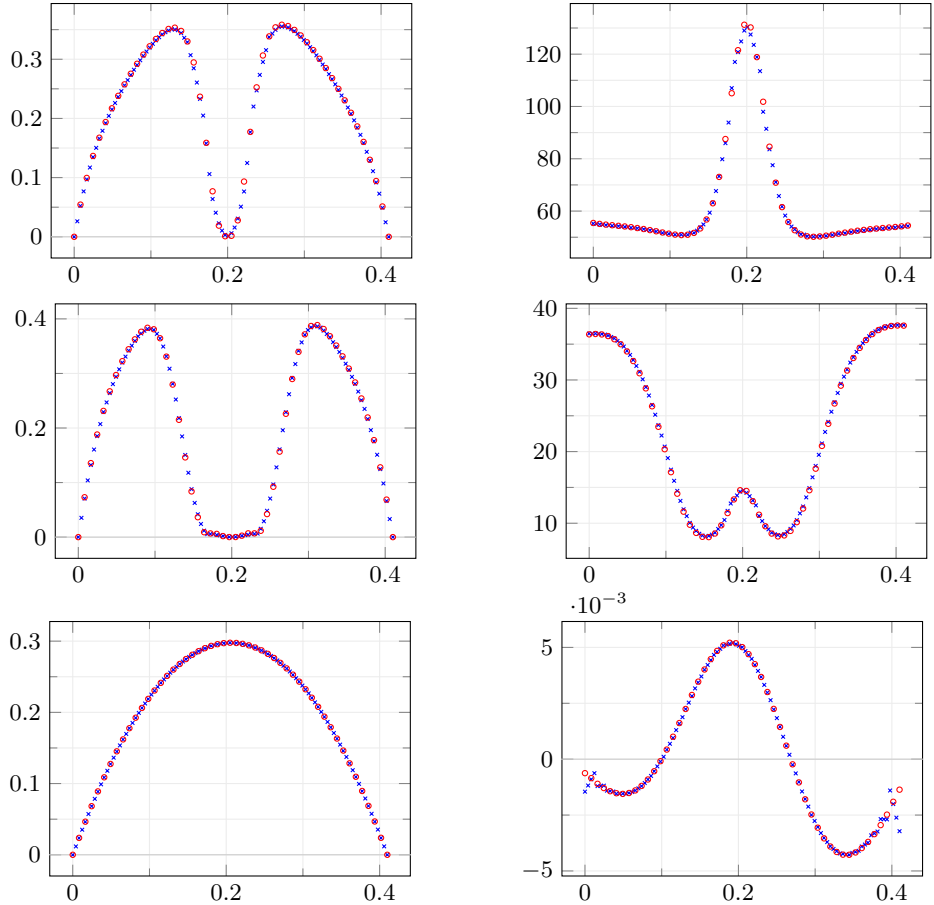

\paragraph{Acknowledgement}
  Anne-Kathrin Wenske and Marc C.\ Steinbach gratefully acknowledge
  the financial support from the Deutsche Forschungsgemeinschaft
  (DFG, German Research Foundation) -- SFB1463 -- 434502799.


\medskip

\bibliographystyle{plain}

\end{document}